\documentclass[a4paper,10pt]{article}

\usepackage{pb-diagram}

\usepackage{amsmath}

\usepackage{theorem}

\usepackage{amscd}

\usepackage{amssymb}

\usepackage{eucal}

\usepackage[T1]{fontenc}

\usepackage[latin1]{inputenc}

\newtheorem{thm}{Theorem}[section]

\newtheorem{prop}{Proposition}[section]

\newtheorem{rem}{Remark}[section]

\theoremstyle{break}

\newcommand{\PP}{\mathbb P}

\newcommand{\XX}{\mathcal{X}}

\newcommand{\hti}{\tilde{h}}

\newcommand{\curv}[2]{\Theta_{#1}(#2)}

\newcommand{\To}{\longrightarrow}

\newcommand{\norm}[1]{\left\Vert#1\right\Vert}

\title{Invariance for multiples of the twisted canonical bundle}
\author{Benoît Claudon}
\date{}

\begin{document}

\maketitle

\section*{Introduction}
Let us consider the following situation : let $\pi:\XX\To \Delta$ a smooth projective family (over the unit disc) et let $L\To \XX$ be a line bundle over $\XX$ endowed with a (possibly singular \footnote[1]{the local weights $\varphi$ of the metric satisfy $\varphi\in L^1_{loc}$ so that $\curv{\hti}{L}=i\partial\overline{\partial}\varphi$ is well defined as a current}) hermitian metric $\hti$ such that :

\begin{description}
\item[(i)]  $\curv{\hti}{L}\geq 0$ as a current ($i.e.$ $(L,\hti)$ is pseudo-effective).
\item[(ii)] the restriction $\hti_{\XX_0}$ of $\hti$ to the central fiber $\XX_0$ is well defined ($i.e.$ if $\varphi$ is a local weight of the metric $\hti$, $\varphi_{\vert \XX_0}\not\equiv -\infty$ and $\varphi_{\vert \XX_0}\in L^1_{loc}$).
\item[(iii)] the multiplier ideal sheaf $\mathcal{I}(\XX_0,\hti_{\XX_0})$ is trivial : $\mathcal{I}(\XX_0,\hti_{\XX_0})=\mathcal{O}_{\XX_0}$.
\end{description}

In this paper, we establish the proof of the following extension result :

\begin{thm}\label{th1}
let $\XX\To \Delta$ a smooth projective family, $m\geq 1$ an integer and let $(L,\hti)$ a hermitian line bundle satisfying the conditions $(i)$, $(ii)$ and $(iii)$ above. Then every section of $m(K_{\XX_0}+L)$ (over $\XX_0$) extends to $\XX$ ; in other words, the restriction map :
$$H^0(\XX,m(K_{\XX}+L))\To H^0(\XX_0,m(K_{\XX_0}+L))$$
is surjective.
\end{thm}

\noindent This result is a "family" version of a result obtained by S. Takayama (see \cite[th 4.1]{T05}).

\vspace{1\baselineskip}

The strategy employed to prove theorem \ref{th1} is the one given by M. Paun to simplify Siu's proof of the invariance of plurigenera and, in the same time, to improve this result. Indeed, in \cite{P05}, M. Paun was able to replace the $L^{\infty}$ hypothesis (originally formulated by Siu to solve the invariance of plurigenera) by an $L^2$ one, which is closely related to extension problems by the way of the Ohsawa-Takegoshi theorem (theorem \ref{th2} below). We would like to point out here the main steps of this method.

First we fix $s\in H^0(\XX_0,m(K_{\XX_0}+L))$, the section we want to extend and let $A$ be an ample line bundle over $\XX$ which satisfies the following conditions :

\begin{description}
\item[(1)] $p(K_{\XX}+L)+A$ is generated by its global sections, say $(s^{(p)}_j)_{j=1..N_p}$, for $0\leq p\leq m-1$
\item[(2)] every section of $m(K_{\XX_0}+L)+A$ extends to $\XX$.
\end{description}
Such an $A$ exists as (1) is required only for finitely many line bundles, and (2) is a consequence of the vanishing of the relevant 
$H^1$ group by Serre's theorem. What we have to do next is to extend the following sections : $s^{\otimes k}\otimes s^{(p)}_j$ with $k \geq 1$, $0\leq p\leq m-1$ and $0\leq j \leq N_p$ to obtain $\widetilde{s}^{(km+p)}_j$, sections of \\
$(km+p)(K_{\XX}+L)+A$. Here the heart of the proof rests on the Ohsawa-Takegoshi theorem : indeed, it implies that we can find such extension with uniform $L^2$ estimates. Then, we use the family $(\widetilde{s}^{(km)}_j)_{j=1..N_p}$ to construct some metrics on $km(K_{\XX}+L)+A$ and the final step is to extract roots and pass to the limit ($i.e.$ dividing by $k$, we consider $m(K_{\XX}+L)+\frac{1}{k}A$ and, passing to the limit, we produce a final metric $h_\infty$ on $m(K_{\XX}+L)$. The main point is that the $L^2$ estimates produce effective bounds and that is why we can pass to the limit. At the end of the proof, the metric $h_\infty$ is used to apply (for the last time !) the Ohsawa-Takegoshi theorem to extend the section $s$.

\vspace{1\baselineskip}

\noindent\textbf{Aknowledgement : } I am very grateful to Mihai Paun for explaining to me his very beautiful method and encouraging me to write down the proof of theorem \ref{th1}.For this and also for many interesting discussions, I would like to thank him.

\section{Preliminaries}

We recall here some facts we need in the proof of the theorem \ref{th1} and fix some notations.

To start with, we want to recall how to define a metric on a line bundle $E\To \XX$ when a family $(s_j)$ of sections of $E$ is given : fix any hermitian (smooth) metric $h$ on $E$ and, for $\sigma \in E$, define :
$$\norm{\sigma}^2=\frac{\norm{\sigma}_h^2}{\sum_j \norm{s_j}_h^2}$$
This (singular) metric is clearly independent of $h$ and its singularities are concentrated along the common zeroes of the sections $(s_j)$ ; moreover, the curvature current of this metric is a closed positive current.

\vspace{1\baselineskip}

As noticed in the introduction above, the main tool of the proof is the $L^2$-extension theorem of Ohsawa and Takegoshi (see \cite{OT87}). However, the version used in the sequel is the one established by Y.-T. Siu in \cite{S02} :

\begin{thm}[Ohsawa-Takegoshi, Siu]\label{th2}
Let $\XX\To \Delta$ a smooth projective family and $L\To \XX$ a line bundle endowed with a (possibly singular) metric $h$ with semi-positive curvature current. Then there exists a (universal) constant $C_0$ such that for every section $\sigma_0\in H^0(\XX_0,K_{\XX_0}+L)$ satisfying :
$$\int_{\XX_0}\norm{\sigma}_h^2<+\infty\, ,$$
there exists $\tilde{\sigma}\in H^0(\XX,K_{\XX}+L)$ with $\tilde{\sigma}_{\vert \XX_0}=\sigma\wedge dt$ and moreover :
$$\int_{\XX}\norm{\tilde{\sigma}}_h^2\leq C_0 \int_{\XX_0}\norm{\sigma}_h^2$$
\end{thm}
The version established in \cite{S02} is actually more general, but the previous statement is enough for our purpose. The crucial point is that the constant $C_0$ is universal : it is independent of $(L,h)$ (for a precise value of $C_0$ see \cite[th. 3.1, p. 241]{S02}).

\vspace{1\baselineskip}

We fix some more notations : we use the ones in the introduction for $s\in H^0(\XX_0,m(K_{\XX_0}+L))$, for $A$ and for the sections $s^{(p)}_{j}\in H^0(\XX,p(K_{\XX}+L)+A)$. If $\omega$ is a hermitian metric on $\XX$, $h_\omega$ will denote the metric induced by $\omega$ on $K_\XX$. Let $h$ a smooth metric on $L$ and $h_A$ a smooth metric on $A$ with $\curv{h_A}{A}>0$ ; if $q\geq 1$ is an integer, $h_q$ will denote the metric $(h_{\omega}\otimes h)^{\otimes q}\otimes h_A$ on $q(K_\XX+L)+A$ (when needed, $h_{q,r}$ will denote the metric $h_{\omega}^{\otimes q}\otimes h^{\otimes r}\otimes h_A$).\\
Consider the metric $\hti$ on $L$ : we can write $\hti=e^{-\tilde{\varphi}}h$ and the assumption on the curvature of $(L,\hti)$ is 
$$\curv{\hti}{L}=\curv{h}{L}+i\partial\overline{\partial}\tilde{\varphi}\geq 0$$
as currents on $\XX$. In particular, this implies that the weight function $\tilde{\varphi}$ is locally bounded from above.

\begin{rem}\label{r1}
the hypothesis made on $\mathcal{I}(\XX_0,\hti_{\XX_0})$ (its triviality) can be expressed in the following way : 
$$\int_{\XX_0}e^{-2\tilde{\varphi}}dV_\omega<+\infty$$
We will denote by $C_L$ this constant in the sequel.
\end{rem}

\section{Proof of the theorem}

As we pointed out in the introduction, we will need precise $L^2$ estimates to achieve passing to the limit ; actually, theorem \ref{th1} will be a straightforward consequence of the following proposition :

\begin{prop}\label{p1}
There exists a constant $C>0$ such that, for all $k\geq 1$, $0\leq p \leq m-1$ and $0\leq j \leq N_p$, there exist some sections $$\widetilde{s}^{(km+p)}_j \in H^0(\XX, (km+p)(K_\XX+L)+A)$$
with $\widetilde{s}^{(km+p)}_{j\vert \XX_0}=s^k\otimes s^{(p)}_j$ and with the following estimates :
\begin{description}
\item[(E1)] if $1\leq p\leq m-1$, we have
$$\int_\XX \frac{\sum_{j=1}^{N_p}\norm{\widetilde{s}^{(km+p)}_j}_{h_{km+p}}^2}{\sum_{j=1}^{N_p-1}\norm{\widetilde{s}^{(km+p-1)}_j}_{h_{km+p-1}}^2}dV_\omega \leq C$$
\item[(E2)] for $p=0$ (and $k\geq 2$), the estimate becomes
$$\int_\XX \frac{\sum_{j=1}^{N_0}\norm{\widetilde{s}^{(km)}_j}_{h_{km}}^2}{\sum_{j=1}^{N_m-1}\norm{\widetilde{s}^{((k-1)m+m-1)}_j}_{h_{(k-1)m+m-1}}^2}dV_\omega \leq C$$
\end{description}
\end{prop}

\noindent\textbf{\large Proof of the proposition \ref{p1} :}\\
To start with, we can consider the sections $s\otimes s^{(0)}_j$ $(0\leq j\leq N_0)$ ; using the poperty $(2)$ of $A$, each of the previous sections extends over $\XX$. Thus, we get the extensions $\widetilde{s}^{(m)}_j $.

Before going further in the proof, it can be useful to do the following remark : by the global property $(1)$ of $A$ (and possibly shrinking $\Delta$), there exists a constant $C_1$ such that
\begin{equation}\label{glob}
\underset{r,q}{\textrm{max}}\,\underset{\XX}{\textrm{sup}}\Big(\frac{\sum_{j=1}^{N_r}\norm{s^{(r)}_j}_{h_r}^2}{\sum_{j=1}^{N_q}\norm{s^{(q)}_j}_{h_q}^2}\Big)\leq C_1
\end{equation}
To prove proposition \ref{p1}, we will proceed inductively and construct the desired extensions step by step ; to this end, we consider the following constant :
$$\widetilde{C}=\textrm{max}(1,\norm{s}^2_{L^\infty,(h_{\omega}\otimes h)^{\otimes m}})C_0C_1C_Le^{2M}
$$
where $M$ is an upper bound for $\tilde{\varphi}$ (we already shrinked $\Delta$ so $M$ exists) and $$\norm{s}_{L^\infty,(h_{\omega}\otimes h)^{\otimes m}}=\underset{x\in\XX_0}{\textrm{sup}}(\norm{s(x)}_{(h_{\omega}\otimes h)^{\otimes m}})$$
We can now initiate the inductive process : to get the extension of the sections $s\otimes s^{(1)}_j$, we consider the line bundle $m(K_\XX+L)+A+L$ we endowed with the metric defined by the family $(\widetilde{s}^{(m)}_j)_{j=0..N_0}$ twisted with the metric $\hti$. This metric has clearly a semi-positive curvature current and, using (\ref{glob}), we have
\begin{equation}\label{eq1}
\frac{\norm{s\otimes s^{(1)}_j}_{h_{m+1,m}\otimes \hti}^2}{\sum_{q=0}^{N_0}\norm{s\otimes s^{(0)}_q}_{h_m}^2}=
\frac{\norm{s\otimes s^{(1)}_j}_{h_{m+1}}^2}{\sum_{q=0}^{N_0}\norm{s\otimes s^{(0)}_q}_{h_m}^2}e^{-2\tilde{\varphi}}\leq C_1 e^{-2\tilde{\varphi}}
\end{equation}
Integrating (\ref{eq1}) over $\XX_0$ and using the remark \ref{r1}, we get
\begin{equation}\label{eq2}
\int_{\XX_0}\frac{\norm{s\otimes s^{(1)}_j}_{h_{m+1,m}\otimes \hti}^2}{\sum_{q=0}^{N_0}\norm{s\otimes s^{(0)}_q}_{h_m}^2}dV_\omega \leq C_1 C_L <+\infty
\end{equation}
We can thus apply the theorem \ref{th2} and we get $\widetilde{s}^{(m+1)}_j$ an extension of $s\otimes s^{(1)}_j$ with the estimate :
\begin{equation}\label{eq3}
\int_{\XX}\frac{\norm{\widetilde{s}^{(m+1)}_j}_{h_{m+1,m}\otimes \hti}^2}{\sum_{q=0}^{N_0}\norm{\widetilde{s}^{(m)}_q}_{h_m}^2}dV_\omega\leq C_0 C_1 C_L
\end{equation}
To have an estimate involving only the metric $h_{m+1}$, we just have to remember that the function $\tilde{\varphi}$ is bounded from above by $M$, so that :
\begin{equation}\label{eq4}
\int_{\XX}\frac{\norm{\widetilde{s}^{(m+1)}_j}_{h_{m+1}}^2}{\sum_{q=0}^{N_0}\norm{\widetilde{s}^{(m)}_q}_{h_m}^2}dV_\omega\leq C_0 C_1 C_Le^{2M}\leq \widetilde{C}
\end{equation}

\noindent Suppose we have already constructed the extension $\widetilde{s}^{(km+p)}$ (with $(k,p)\neq (1,0)$) with the desired estimates ; we now have to climb to the next step. To do this, we separate the two different following case :

\vspace{1\baselineskip}

\noindent\textbf{case 1 : $p<m-1$}\\
we consider the line bundle $(km+p)(K_\XX+L)+A+L$ that we endowed with the metric coming from the family $(\widetilde{s}^{(km+p)}_q)_{q=0\dots N_{p}}$ twisted by $\hti$ ; as in the case treated above, we have the following estimates on $\XX_0$ :
\begin{equation}
\frac{\norm{s^k\otimes s^{(p+1)}_j}_{h_{km+p+1,km+p}\otimes \hti}^2}{\sum_{q=0}^{N_p}\norm{ \widetilde{s}^{(km+p)}_q}_{h_{km+p}}^2}=
\frac{\norm{s^k\otimes s^{(p+1)}_j}_{h_{km+p+1}}^2}{\sum_{q=0}^{N_p}\norm{s^k\otimes s^{(p)}_q}_{h_{km+p}}^2}e^{-2\tilde{\varphi}}\leq C_1 e^{-2\tilde{\varphi}}
\end{equation}
and we can then extend $s^k\otimes s^{(p+1)}_j$ with estimate, exactly in the same way as in the first step of the induction.

\vspace{1\baselineskip}

\noindent\textbf{case 2 : $p=m-1$}\\
we still have to consider the line bundle $(km+m-1)(K_\XX+L)+A+L$ endowed with the metric coming from the family $(\widetilde{s}^{(km+m-1)}_q)_{q=0\dots N_{m-1}}$ twisted by $\hti$ ; at this step, we obtain the needed estimate (on $\XX_0$) as follows :
\begin{eqnarray*}
\int_{\XX_0}\frac{\norm{s^{k+1}\otimes s^{(0)}_j}_{h_{(k+1)m,km+m-1}\otimes \hti}^2}{\sum_{q=0}^{N_{m-1}}\norm{\widetilde{s}^{(km+m-1)}_q}_{h_{km+m-1}}^2}dV_\omega & = & \int_{\XX_0}\frac{\norm{s^{k+1}\otimes s^{(0)}_j}_{h_{(k+1)m}}^2}{\sum_{q=0}^{N_{m-1}}\norm{s^k\otimes s^{(m-1)}_q}_{h_{km+m-1}}^2}e^{-2\tilde{\varphi}}dV_\omega \\
 & \leq & C_1\int_{\XX_0}\norm{s}_{(h_{\omega}\otimes h)^{\otimes m}}^{2}e^{-2\tilde{\varphi}}dV_\omega \\
 & \leq & C_1C_L\norm{s}^2_{L^\infty,(h_{\omega}\otimes h)^{\otimes m}}
\end{eqnarray*}
Applying theorem \ref{th2}, we find a section $\widetilde{s}^{((k+1)m)}_j \in H^0(\XX,(k+1)m(K_\XX+L)+A)$ with $\widetilde{s}^{((k+1)m)}_{j\vert \XX_0}=s^{k+1}\otimes s^{(0)}_j$ and
\begin{equation}
\int_{\XX}\frac{\norm{\widetilde{s}^{((k+1)m)}_j}_{h_{(k+1)m,km+m-1}\otimes \hti}^2}{\sum_{q=0}^{N_{m-1}}\norm{\widetilde{s}^{(km+m-1)}_q}_{h_{km+m-1}}^2}dV_\omega \leq C_0C_1C_L\norm{s}^2_{L^\infty,(h_{\omega}\otimes h)^{\otimes m}}
\end{equation}
In order to get the final inductive estimate, we use again the fact that $\tilde{\varphi}$ is bounded from above by $M$ and then
\begin{equation}\label{eq6}
\int_{\XX}\frac{\norm{\widetilde{s}^{((k+1)m)}_j}_{h_{(k+1)m}}^2}{\sum_{q=0}^{N_{m-1}}\norm{\widetilde{s}^{(km+m-1)}_q}_{h_{km+m-1}}^2}dV_\omega \leq e^{2M}C_0C_1C_L\norm{s}^2_{L^\infty,(h_{\omega}\otimes h)^{\otimes m}}\leq\widetilde{C}
\end{equation}
We just have to pose $C=\widetilde{C}\cdot\textrm{max}(N_0,\dots ,N_{m-1})$ to conclude the proof of proposition \ref{p1}.$\square$

\newpage

\noindent\textbf{\large Proof of theorem \ref{th1} :}\\

The end of the proof is now reduced to extract roots of the metrics induced by the families $(\widetilde{s}^{(km+p)}_q)_{q=0\dots N_p}$ (see also \cite{P05}) ; indeed, we consider the following weight functions :
$$f_k=\frac{1}{2}\log (\sum_{j=1}^{N_0}\norm{\widetilde{s}^{(km)}_j}_{h_{km}}^{2})$$
Possibly shrinking the disk $\Delta$ (to use Jensen inequality and to bound the $L^2$ norms of $\widetilde{s}^{(m)}_j$), the inductive estimates $(E1)$ and $(E2)$ in the proposition \ref{p1} and the concavity of the logarithm function implies the following inequalities :
\begin{equation}\label{int}
\frac{1}{k}\int_{\XX} f_kdV_\omega \leq C'
\end{equation}
where $C'$ is a positive constant (independent of $k$). Moreover, $f_k$ satisfy the properties :
\begin{equation}\label{curv}
\curv{h_m}{m(K_\XX+L)}+\frac{i}{k}\partial\overline{\partial}f_k\geq -\frac{1}{k}\curv{h_A}{A}
\end{equation}
(in the sense of currents) and, on the central fiber, we have 
\begin{equation}\label{p-est}
\frac{2}{k}f_{k\vert \XX_0}=\log(\norm{s}^2)+\frac{1}{k}\log(\sum_{j=1}^{N_0}\norm{s^{(0)}_j}_{h_0}^2)
\end{equation}
together with the mean value inequality, (\ref{int}) and (\ref{curv}) imply the existence of uniform local upper bounds for the functions $\frac{1}{k}f_k$ (on each relatively compact subset of $\XX$) and thus we can consider :
$$f_\infty=\overline{\underset{k\to +\infty}{\lim \textrm{reg}}}\,\frac{1}{k}f_k$$
the upper semi-continuous enveloppe of the family $(\frac{1}{k}f_k)_{k\geq 1}$ : this is still a quasi-psh function on $\XX$. The property (\ref{p-est}) yields the pointwise estimate (on the central fiber $\XX_0$) :
\begin{equation}\label{p-est2}
\norm{s}^2e^{-2f_\infty}\leq 1
\end{equation}
The metric $h_\infty =e^{-f_\infty}h_m$ is now a (singular) metric with semi-positive current of curvature (by property (\ref{curv}), after passing to the limit) and $s$ is bounded for this metric. To conclude the proof, we consider the metric $g=h_{\infty}^{\frac{m-1}{m}}\otimes \hti$ on the line bundle $(m-1)(K_\XX+L)+L$ ; this is still a metric with semi-positive curvature and the Hölder inequality gives
\begin{eqnarray*}
\int_{\XX_0} \norm{s}_g^2 & = & \int_{\XX_0} \norm{s}^2e^{-2\frac{(m-1)}{m}f_\infty -2\tilde{\varphi}}dV_\omega \\
 & = & \int_{\XX_0} \norm{s}^{2\frac{(m-1)}{m}}e^{-2\frac{(m-1)}{m}(f_\infty +\tilde{\varphi})} \norm{s}^{\frac{2}{m}}e^{-\frac{2}{m}\tilde{\varphi}}dV_\omega \\
  & \leq & \Big( \int_{\XX_0} \norm{s}^2e^{-2f_\infty}e^{-\tilde{2\varphi}}dV_\omega\Big)^{\frac{m-1}{m}}
  \Big(\int_{\XX_0} \norm{s}^2e^{-2\tilde{\varphi}}dV_\omega\Big)^{\frac{1}{m}}
\end{eqnarray*}
Using (\ref{p-est2}) and the remark \ref{r1}, we see that $s$ is actually $L^2$ for the metric $g$. We can thus apply a last time the Ohsawa-Takegoshi theorem \ref{th2} and then obtain the desired extension of $s$.$\square$

\section{Further extension results}

At this stage, we can combine different kinds of extension results to obtain some quite general statements. Let us first recall the following theorem stated by M. Paun in \cite{P05} :
\begin{thm}[Paun]\label{th3}
Let $\XX\To \Delta$ a smooth projective family, $m\geq 1$ an integer and let $(L,h)$ a hermitian line bundle over $\XX$ such that its curvature satisfy : $\curv{h}{L}\geq 0$ (as a current) and such that the restriction of $h$ to the central fiber $\XX_0$ is well defined. Then any section of $(mK_{\XX_0}+L)\otimes\mathcal{I}(h_{\XX_0})$ extends to $\XX$.
\end{thm}
The way of proving this theorem is exactly the same as for theorem \ref{th1} : actually (as already noticed in the introduction), our proof of theorem \ref{th1} is directly inspired from this method.\\

Now, if $(L,h)$ is a (singular) hermitian line bundle over $\XX$ with a semipositive curvature current, the following statement is a kind of interpolation of theorems \ref{th1} and \ref{th3} :

\begin{thm}\label{th4}
Let $\XX\To \Delta$ a smooth projective family, $m,p\geq 1$ integers and let $(L,h)$ a hermitian line bundle over $\XX$ as above such that the restriction of $h$ to the central fiber $\XX_0$ is well defined. Assume moreover that the following condition holds : $\mathcal{I}(h_{\XX_0}^q)=\mathcal{O}_{\XX_0}$ where $p=(m-1)q+r$ (with $0\leq r\leq m-2$). Then, any section of $(mK_{\XX_0}+pL)\otimes\mathcal{I}(h_{\XX_0}^r)$ extends to $\XX$.
\end{thm}
The reason for which we have to write $p=(m-1)q+r$ rather than $p=mq+r$ is the following : the induction process is a sequence of sub-process, each of them divided into $m$ steps. Thus, using the method above, the triviality of $\mathcal{I}(h_{\XX_0}^q)$ allows us to apply Ohsawa-Takegoshi theorem in the first $(m-1)$ steps and, for the final step, the section has to be $L^2$ with respect to $h^r$. That is why we have to consider the decomposition $p=(m-1)q+r$.\\

Actually, as it was pointed out by J.-P. Demailly (\cite{D06}), we can consider mixed problems of extension of pluricanonical sections : 

\begin{thm}[Demailly]\label{th5}
Let $\XX\To \Delta$ a smooth projective family, $m\geq 1$ an integer and let $(L_j,h_j)_{0\leq j\leq m-1}$ be hermitian line bundles over $\XX$ with semipositive curvature current $\curv{h_j}{L_j}\geq 0$. Assume that :
\begin{description}
\item[(i)]  the restriction of $h_j$ to the central fiber $\XX_0$ is well-defined
\item[(ii)] for $1\leq j\leq m-1$, the multiplier ideal sheaf $\mathcal{I}(h_{j\vert\XX_0})$ is trivial
\end{description}
Then, any section of $(mK_{\XX_0}+\sum_j L_j)\otimes\mathcal{I}(h_{0\vert\XX_0})$ over the central fiber of the family extends to $\XX$.
\end{thm}
For instance, theorem \ref{th5} applied to $L_0=rL$ and $L_j=qL$ for $1\leq j\leq m-1$ is nothing but theorem \ref{th4} above.

\newpage

\section{Comparison with the projective case}

As noticed at the beginning of this paper, theorem \ref{th1} is a family version of another result of S. Takayama ; actually, the proof given here can be immediately adapted to obtain the following statement of this result :

\begin{thm}[Takayama]\label{th6}
Let $X$ be a smooth projective manifold, $S\subset X$ a smooth irreducible hypersurface and $L$ a line bundle over $X$ endowed with a singular metric $h$ such that :
\begin{description}
\item[(i)]  $\curv{h}{L}\geq \epsilon\omega$ (with $\epsilon>0$ and $\omega$ any smooth hermitian metric on $X$)
\item[(ii)] the restriction $h_S$ of the metric $h$ to $S$ is well defined and $\mathcal{I}(S,h_S)=\mathcal{O}_{S}$
\end{description}
Then, for any integer $m\geq 1$, the natural restriction map :
$$H^0(X,m(K_X+S+L))\To H^0(S,m(K_S+L))$$
is surjective.
\end{thm}
In this setting, the Ohsawa-Takegoshi theorem \ref{th2} is however no longer valid so, instead of it, we have to use the following extension result (which is a simple consequence of the Nadel vanishing theorem) :

\begin{prop}\label{p2}
Let $X$ be a smooth projective manifold, $S\subset X$ a smooth irreducible hypersurface and $(L,h)$ a singular hermitian line bundle over $X$ satisfying :
\begin{description}
\item[(i)] $\curv{h}{L}\geq \epsilon\omega$ 
\item[(ii)] $h_S$ is well defined. 
\end{description}
Then, for every section $\sigma\in H^0(S,(K_S+L)\otimes\mathcal{I}(h_S))$, there exists a section $\tilde{\sigma}\in H^0(X,K_X+S+L)$ which extends $\sigma$ over $X$.
\end{prop}
Here, we can remark the following : theorem \ref{th2} and proposition \ref{p2} correspond both to the case $m=1$ in the different extension results for pluricanonical forms.

\vspace{0.5cm}

Actually the main difference between theorem \ref{th1} and theorem \ref{th6} sits in the positivity assumption for the line bundle $L$ : in the projective case, we have to require strict positivity for $L$. The reason is the following : as in the family setting, we try to extend some sections $\sigma^k\otimes s^{(p)}_j$ but using proposition \ref{p2} instead of the Ohsawa-Takegoshi theorem \ref{th2} ; thus, we cannot use a limit process to extract roots and the strict positivity of $L$ is essential to balance the negative contribution of $-\frac{1}{k}A$ (where $A$ is the auxiliary ample line bundle and $k$ is chosen big enough). Then, this emphasizes the key role played by the Ohsawa-Takegoshi theorem : extending sections with precise $L^2$ estimates.\\

As a final remark, we can wonder if other (weakened) positivity assumptions on $L$ and $S$ (instead of $(i)$ and $(ii)$ in theorem \ref{th6}) can lead to the same conclusion : for instance, is it true that nefness of $L-S$ implies the surjectivity of the restriction map ? The answer to the preceeding question is actually negative as the following example shows (see also \cite{DPS94}) :\\

Let $E$ be an elliptic curve and $V$ be the rank 2 vector bundle over $E$ defined as the (unique) non split extension :
$$0\To \mathcal{O}_E \To V\To \mathcal{O}_E \To 0$$
In particular, $V$ is numerically flat : $c_1(V)=0$ and $c_2(V)=0$. Now, consider the ruled surface $X=\PP(V)$ and the corresponding section $S=\PP(\mathcal{O}_E)\subset X$. It is an easy matter to check that $S$ satisfy the following :
$$S^2=0\,,\quad \mathcal{O}_X(S)=\mathcal{O}_{\PP(V)}(1)\,,\quad \mathcal{O}_S(S)=\mathcal{O}_S$$
Moreover, the canonical bundle of $X$ is given by :
$$K_X=\mathcal{O}_X(-2S)$$
Now choose $L$ be the line bundle : $L=\mathcal{O}_X(2S)=\mathcal{O}_{\PP(V)}(2)$. $V$ being numerically flat, it is a nef vector bundle and thus 
$$L-S=\mathcal{O}_X(2S)-\mathcal{O}_X(S)=\mathcal{O}_{\PP(V)}(1)$$
is nef too. Furthermore, we have :
\begin{align*}
K_X+L+S &=\mathcal{O}_X(-2S)+\mathcal{O}_X(2S)+\mathcal{O}_X(S)=\mathcal{O}_{\PP(V)}(1)\\
K_S+L_{\vert S} &=(K_X+L+S)_{\vert S}=\mathcal{O}_S(S)=\mathcal{O}_S
\end{align*}
It is now clear that, for $m\geq 1$, the restriction map :
$$H^0(X,\mathcal{O}_{\PP(V)}(m))\simeq H^0(X,m(K_X+S+L))\To H^0(S,m(K_S+L))\simeq H^0(S,\mathcal{O}_S)$$
cannot be surjective.\\

With this example, it should be clear that the relationship between the positivity of $L$ and $S$ play a crucial role in the problem of extending pluricanonical sections from subvariety to the ambiant space.

\begin{rem}
In \cite{DPS94}, the line bundle $L$ is actually an example of nef line bundle which however does not admit any smooth metric with semipositive curvature ($i.e.$ $L$ is nef but not hermitian semipositive).
\end{rem}

\bibliographystyle{amsalpha}
\bibliography{biblio}

\vspace*{0.5cm}
\begin{flushright}
\begin{minipage}{5cm}
Beno\^it CLAUDON\\
Universit\'e Nancy 1\\
Institut Elie Cartan\\
BP 239\\
54 506 Vandoeuvre-l\`es-Nancy\\
Cedex (France)

\vspace*{0.3cm}
\noindent Benoit.Claudon@iecn.u-nancy.fr
\end{minipage}
\end{flushright}

\end{document}